\documentclass[11pt,reqno]{article}

\usepackage{amsfonts}
\usepackage{amsmath}
\usepackage{amsthm}
\usepackage{amssymb}
\usepackage{mathrsfs}
\usepackage{amstext}
\usepackage{graphicx}
\font\msbm=msbm10

\numberwithin{equation}{section}

\theoremstyle{plain}
\newtheorem{Theorem}{Theorem}[section]
\newtheorem{lemma}[Theorem]{Lemma}

\newtheorem{proposition}[Theorem]{Proposition}

\newtheorem{remark}[Theorem]{Remark}
\def\mathbb#1{\hbox{\msbm{#1}}}

\newcommand{\N}{{\mathbb{N}}}
\newcommand{\R}{{\mathbb{R}}}

\newcommand{\beq}{\begin{eqnarray}}
\newcommand{\eeq}{\end{eqnarray}}
\newcommand{\beqn}{\begin{eqnarray*}}
\newcommand{\eeqn}{\end{eqnarray*}}

\newcommand{\supp}{\operatorname{supp}}

\renewcommand{\qed}{\rule{2.5mm}{2.5mm}}
\newenvironment{Proof}{\noindent
{\bf\underline{Proof:} }}
{\hspace*{\fill}\qed\vskip1em}
\begin{document}
\title{The Gelfand widths of $\ell_p$-balls for $0<p\leq 1$}

\author{Simon Foucart, Alain Pajor, Holger Rauhut, Tino Ullrich}

\date{April 19 2010}
\maketitle

\begin{abstract} We provide sharp lower and upper bounds for the Gelfand widths 
of $\ell_p$-balls in the $N$-dimensional $\ell_q^N$-space 
for $0<p\leq 1$ and $p<q \leq 2$. Such estimates 
are highly relevant to the novel theory of compressive sensing, and our proofs rely on
methods from this  area.  
\end{abstract}

\begin{tabbing}
{\bf Key Words:} Gelfand widths, compressive sensing, sparse recovery, 
$\ell_1$-minimization,\\ $\ell_p$-minimization.
\end{tabbing}
{\bf AMS Subject classification:} 41A46, 46B09.

\section{Introduction}

Gelfand widths are an important concept in classical and modern approximation and complexity theory. 
They have found 
recent interest in the 
rapidly emerging field of compressive sensing \cite{carota06,codade09,do06-2} 
because they give general performance bounds 
for sparse recovery methods. Since vectors in $\ell_p$-balls, $0<p\leq 1$, can be well-approximated by sparse
vectors, the Gelfand widths of such balls are particularly relevant in this context.
In remarkable papers \cite{ka77,gl84,gagl84} from the 1970s and 80s due to Kashin, Gluskin, and Garnaev, upper and lower estimates for the Gelfand widths of $\ell_1$-balls are provided. 
In his seminal paper introducing compressive sensing \cite{do06-2}, Donoho extends these estimates
to the Gelfand widths of $\ell_p$-balls with $p<1$. Unfortunately, his proof of the lower bound contains a gap.
In this article, we address 
this issue by supplying a complete proof. 
To this end, we  proceed in an entirely different way than Donoho. 
Indeed, we use compressive sensing methods to establish the lower bound in a more 
intuitive way. Our method is new even for the case $p=1$.
For completeness, we also
give a proof of the upper bound based again on compressive sensing arguments. These arguments also provide the same sharp asymptotic behavior for the Gelfand widths of weak-$\ell_p$-balls.

\subsection{Main Result}

In this paper, we consider the finite-dimensional  spaces $\ell_p^N$, that is, $\R^N$ endowed with the usual $\ell_p$-(quasi-)norm defined, for  $x \in \R^N$, by 
\[
\|x\|_p := \Big( \sum_{\ell=1}^N |x_\ell|^p\Big)^{1/p}, \qquad 0<p<\infty,
\qquad \quad
\|x\|_\infty := \max_{\ell=1,\hdots,N} |x_\ell|.
\] 
For $1\leq p \leq \infty$, this is a norm, while for $0<p<1$, it only 
satisfies the $p$-triangle inequality
\begin{equation}\label{ptriangle}
\|x+y\|_p^p \leq \|x\|_p^p + \|y\|_p^p,
\qquad x,y \in \R^N.
\end{equation}
Thus, $\|\cdot\|_p$ defines a quasi-norm with 
(optimal) quasi-norm constant $C = \max\{1,2^{1/p-1}\}$.
As a reminder, $\|\cdot\|_X$ is called a quasi-norm on $\R^N$ with quasi-norm constant $C\geq 1$ if it obeys 
the quasi-triangle inequality 
$$
   \|x+y\|_X \leq C(\|x\|_X + \|y\|_X),
   \qquad x,y \in \R^N.
$$ 
Other quasi-normed spaces considered in this paper are the spaces weak-$\ell_p^N$,
that is, $\R^N$ endowed with the $\ell_{p,\infty}$-quasi-norm defined, for  $x \in \R^N$, by
$$
      \|x\|_{p,\infty} := \max\limits_{\ell = 1,\hdots,N} \ell^{1/p}|x_{\ell}^{\ast}|,
      \qquad 0<p \le \infty,
$$
where $x^{\ast}\in \R^N$ is a non-increasing rearrangement of $x$.  
 We shall investigate the Gelfand widths in  $\ell_q^N$ of the unit balls
$B_p^N := \{x \in \R^N, \|x\|_p \leq 1\}$ and
$B_{p,\infty}^N := \{x \in \R^N, \|x\|_{p,\infty} \leq 1\}$ of $\ell_p^N$ and $\ell_{p,\infty}^N$  for $0 < p \le 1$ and $p < q \le 2$.

We recall that the Gelfand width of order $m$
of a subset $K$ of $\R^N$ 
in the (quasi-)normed space $(\R^N, \|\cdot\|_X)$ is defined as
\[
d^m(K,X) := \inf_{A \in \R^{m \times N}} \sup_{v \in K \cap \ker A} \|v\|_X,
\]
where $\ker A := \{v \in \R^N, Av = 0\}$ denotes the kernel of $A$. 
It is well-known that the above infimum 
is actually realized \cite{pi85}. 
Let us observe that $d^m(K,X) = 0$ for $m\geq N$
when $0 \in K$, so that we restrict our considerations to the case $m<N$ in the sequel.
Let us also observe that the simple inclusion $B_p^N \subset B_{p,\infty}^N$ implies 
$$
d^m(B_p^N,\ell_q^N) \le d^m(B_{p,\infty}^N,\ell_q^N).
$$
From this point on, we aim at finding a lower bound for $d^m(B_p^N,\ell_q^N)$ and an upper bound for $d^m(B_{p,\infty}^N,\ell_q^N)$ with the same asymptotic behaviors. 
Our main result is summarized below.

\begin{Theorem}\label{thm:main} 
For $0<p\leq 1$ and $p < q \leq 2$,  there exist constants 
$c_{p,q}, C_{p,q} > 0$ depending only on $p$ and $q$ such that, if $m < N$, then
\begin{equation}\label{upper:lower:Gelfand}
c_{p,q} \min\left\{1,\frac{\ln(N/m)\hspace{-.4mm}+\hspace{-.4mm}1}{m} \right\}^{1/p-1/q}\hspace{-.5mm} 
\leq \, d^m(B_{p}^N,\ell_q^N) \,
\leq  C_{p,q} \min\left\{1,\frac{\ln(N/m)\hspace{-.4mm}+\hspace{-.4mm}1}{m} \right\}^{1/p-1/q}
\hspace{-.5mm}, 
\end{equation}
and, if $p<1$,
\begin{equation}
\label{upper:lower:Gelfand2}
c_{p,q} \min\left\{1,\frac{\ln(N/m)\hspace{-.4mm}+\hspace{-.4mm}1}{m} \right\}^{1/p-1/q}\hspace{-.5mm}
\leq d^m(B_{p,\infty}^N,\ell_q^N)
\leq  C_{p,q} \min\left\{1,\frac{\ln(N/m)\hspace{-.4mm}+\hspace{-.4mm}1}{m} \right\}^{1/p-1/q}\hspace{-.5mm}.
\end{equation}
\end{Theorem}

In the case $p=1$ and $q=2$, the upper bound of \eqref{upper:lower:Gelfand} with  a slightly worse $\log$-term was shown 
by Kashin in \cite{ka77} by considering Kolmogorov widths, which are dual to Gelfand widths \cite{lomavo96,pi85}. 
The lower bound and the optimal $\log$-term for $p=1$ and $1<q\leq 2$
were provided by Garnaev and Gluskin in \cite{gl84,gagl84}, again via Kolmogorov widths. 
An alternative proof of the upper and lower estimates of \eqref{upper:lower:Gelfand} with $p=1$ 
was given by Carl and Pajor in \cite{capa88}. 
They did not pass to Kolmogorov widths, but rather used Carl's theorem
\cite{ca81-1} (see also \cite{cast90,pi85}) that bounds in particular Gelfand numbers from below by entropy numbers, 
which are completely understood even for $p,q<1$, see 
\cite{sc84-1, guli00, ku01-1}. An upper bound
for $p<1$ and $q=2$ was first provided by Donoho \cite{do06-2} with $\log(N)$ instead of $\log(N/m)$. 
With an adaptation of a method from \cite{lomavo96}, Vyb{\'i}ral \cite[Lem. 4.11]{vy08} also provided the upper 
bound of \eqref{upper:lower:Gelfand} when $0 < p \leq 1$. 
In Section 3, we use  compressive sensing techniques to give an alternative proof that provides the upper bound of \eqref{upper:lower:Gelfand2}.

Donoho's attempt to prove the lower bound of \eqref{upper:lower:Gelfand} for the case 
$0<p<1$ and $q=2$ consists in
applying Carl's theorem and then using known estimates for entropy numbers, 
similarly to the approach by Carl and Pajor for $p=1$. 
However, it is unknown whether Carl's theorem extends to 
quasi-norm balls, in particular to $\ell_p$-balls with $p<1$. The standard proof of Carl's theorem 
for Gelfand widths \cite{cast90, lomavo96}  uses duality arguments, which are not available for quasi-Banach spaces. We believe that Carl's theorem 
actually fails for Gelfand widths of general quasi-norm balls, 
although it turns out to be a posteriori  true  in our specific situation 
due to Theorem \ref{thm:main}. 

We shortly comment on the case $q > 2$. Since then $\|v\|_q \leq \|v\|_2$ for all $v \in \R^N$, we have the upper estimate
\begin{equation}\label{upperq2}
d^m(B_p^N,\ell_q^N) \leq d^m(B_p^N, \ell_2^N)
\leq  C_{p,2} \min\left\{1,\frac{\ln(N/m)+1}{m} \right\}^{1/p-1/2}.
\end{equation}
The lower bound in \eqref{upper:lower:Gelfand} extends to $q>2$, but is unlikely to be optimal in this case. It seems rather
that \eqref{upperq2} is close to the correct behavior. At least for $p=1$ and $q> 2$, \cite{gl82-1} gives lower estimates of
related Kolmogorov widths which then leads to (see also \cite{vy08})
\[
d^m(B_1^N,\ell_q^N) \geq c_q m^{-1/2}\,.
\]
The latter matches \eqref{upperq2} up to the $\log$-factor. We expect a similar behavior for $p<1$, 
but this fact remains to be proven.

\subsection{Relation to Compressive Sensing}

Let us now outline the connection to compressive sensing. This emerging theory explores the
recovery of  vectors $x \in \R^N$ from incomplete linear information $y = Ax \in \R^m$, where
$A \in \R^{m \times N}$ and $m < N$. Without additional information, reconstruction of $x$ from $y$
is clearly impossible since, even in the full rank case, the system $y = Ax$ has infinitely many solutions. 
Compressive sensing makes the additional assumption that $x$ is sparse or at least compressible. A vector $x \in \R^N$ 
is called $s$-sparse if at most $s$ of its coordinates are non-zero. The error of best $s$-term approximation
is defined as 
$$
\sigma_s(x)_p := \inf\{ \|x - z \|_p, z \mbox{ is } s\mbox{-sparse} \}.
$$
Informally, a vector $x$ is called compressible if $\sigma_s(x)_p$ decays quickly in $s$. 
It is classical to show that, for $q > p$,
\begin{align}
\label{compressible}
\sigma_s(x)_q &\leq \frac{1}{s^{1/p-1/q}} \, \|x\|_p\,,\\
\label{compressible2}
 \sigma_s(x)_q & \leq \frac{D_{p,q}}{s^{1/p-1/q}}\, \|x\|_{p,\infty} \,,\qquad \qquad D_{p,q} := (q/p-1)^{-1/q}\,.
\end{align}
This implies that the balls $B_p^N$ and $B_{p,\infty}^N$ with $p\leq 1$ serve as good models for compressible signals: 
the smaller $p$, the better $x$ in $B_p^N$ or in $B_{p,\infty}^N$ is approximable in $\ell_q^N$ by $s$-sparse vectors.

The aim of compressive sensing is to find good pairs of linear measurement maps 
$A \in \R^{m \times N}$ and
(non-linear) reconstruction maps $\Delta : \R^m \to \R^N$
that recover compressible vectors $x$ with small errors $x-\Delta(Ax)$. 
In order to measure the performance of a pair $(A,\Delta)$, one defines, 
for a subset $K$ of $\R^N$ and a (quasi-)norm $\|\cdot\|_X$ on $\R^N$, 
\[
E_m(K,X) := \inf_{A \in \R^{m \times N}, \, \Delta: \R^m \to \R^N} \;\; \sup_{x \in K} \|x - \Delta(Ax)\|_X.
\]
Quantities of this type play a crucial role in the modern field of information based complexity \cite{NoWo08}. 
In our situation, only linear information is allowed in order to recover $K$ uniformly. 
The quantities $E_m(K,X)$ are closely linked to the Gelfand widths, as stated in the following 
proposition \cite{do06-2,codade09},
see also \cite{pi86,no95-1}. 

\begin{proposition}\label{prop:Em} 
Let $K \subset \R^N$ be such that $K=-K$ and $K + K \subset C_1 K$ for some $C_1 \geq 2$,
and let $\|\cdot\|_X$ be a quasi-norm on $\R^N$ with quasi-norm constant $C_2$.
Note that $C_1 = 2$ if $K$ is a norm ball and that $C_2=1$ if $\|\cdot\|_X$ is a norm.
Then
\[
C_2^{-1} d^m(K,X) \leq E_m(K,X) \leq C_1 d^m(K,X). 
\]
\end{proposition}
Combining the previous proposition with Theorem \ref{thm:main} 
gives optimal performance bounds for
recovery of compressible vectors in $B_p^N$, $0<p\leq 1$, when the error is measured in 
$\ell_q$, $p<q\leq 2$. 
Typically, the  most interesting case is $q=2$, for which we end up with
\[
c_p \min\left\{1,\frac{\ln(N/m) + 1}{m} \right\}^{1/p-1/2} \leq E_m(B_p^N, \ell_2^N) \leq C_p \min\left\{1,\frac{\ln(N/m)+1}{m}\right\}^{1/p-1/2}. 
\] 
For practical purposes, it is of course desirable to find matrices $A \in \R^{m \times N}$ and efficiently implementable reconstruction maps $\Delta$ that realize the optimal 
bound above. 
For instance,
Gaussian random matrices $A \in \R^{m \times N}$, i.e.,
matrices whose entries are
independent copies of a zero-mean Gaussian variable,
provide optimal measurement maps with high probability \cite{cata06,do06-2,badadewa08}. 
An optimal reconstruction map is obtained via basis pursuit \cite{chdosa99,do06-2,cata06}, i.e.,  via
the $\ell_1$-minimization mapping given by
\[
\Delta_1(y) := {\rm{arg}} \min \|z\|_1 \quad \mbox{ subject to } Az = y.
\]
This mapping can be computed with efficient convex optimization methods \cite{bova04}, and works very well in practice. The proof of the lower bound in \eqref{upper:lower:Gelfand}
will further involve $\ell_p$-minimization for $0<p\le1$
via the mapping
\[
\Delta_p(y) :={\rm{arg}} \min \|z\|_p \quad \mbox{ subject to } Az = y. 
\]
A key concept in the analysis of sparse recovery via 
$\ell_p$-minimization
is the  restricted isometry property (RIP).
This well-established concept in compressive sensing \cite{cata06,carota06-1}
is the main tool for the proof of the upper bound in \eqref{upper:lower:Gelfand2}.
We recall that the $s$th order restricted isometry constant $\delta_s(A)$ of a matrix $A \in \R^{m \times N}$ is defined as the smallest $\delta>0$ such that
$$
(1-\delta) \|x\|_2^2 \le \|Ax\|_2^2 \le (1+\delta) \|x\|_2^2
\qquad \mbox{for all $s$-sparse } x \in \R^N.
$$
Small restricted isometry constants imply stable recovery by $\ell_1$-minimization, as well as by 
$\ell_p$-minimization for $0<p < 1$.
For later reference, we state the following result \cite{carota06-1,ca08,fola09}.
\begin{Theorem}\label{thm:l1} Let $0<p\leq 1$. 
If $A \in \R^{m \times N}$ has a restricted isometry constant 
$\delta_{2s} < \sqrt{2}-1$,
then, for all $x \in \R^N$,
\begin{equation}
\label{StableRec}
\| x-\Delta_p(Ax)\|_p^p \leq C\sigma_s(x)_p^p,
\end{equation}
where $C>0$ is a constant that depends only on $\delta_{2s}$. 
In particular, 
the reconstruction of $s$-sparse vectors is exact.
\end{Theorem}
Given a prescribed $0 < \delta < 1$,
it is known \cite{badadewa08,cata06,mepato09} 
that, if
the entries of the matrix $A$ are independent copies of a zero-mean Gaussian variable with variance $1/m$, 
then  there exist constants $C_1, C_2 > 0$ (depending only on $\delta$) such that $\delta_s(A) \le \delta$ holds with
probability greater than $1-e^{-C_2 m}$
provided that
\begin{equation}\label{mRIP}
m \ge C_1 s \ln(eN/s).
\end{equation}
In particular, there exists a matrix $A \in \R^{m \times N}$ such that the pair $(A,\Delta_1)$, and more generally $(A,\Delta_p)$ for $0<p\leq 1$, allows stable recovery in the sense of \eqref{StableRec} as soon as the number of measurements satisfies \eqref{mRIP}.
Vice versa, we will see in Theorem \ref{stable}
that the existence of any pair $(A,\Delta)$ allowing such a stable recovery 
forces the number of measurements to satisfy \eqref{mRIP}. 

Lemma \ref{pajor2}, which is of independent interest, estimates the minimal number of measurements for the pair $(A,\Delta_p)$ to allow exact (but not necessarily stable) recovery of sparse vectors.
Namely, we must have 
\begin{equation}
\label{Improve3s}
m \geq c_1ps\ln(N/(c_2 s))
\end{equation}
for some explicitly given constants $c_1,c_2 > 0$. 
In the case $p=1$, this result can be also obtained as a consequence of a corresponding lower bound
on neighborliness of centrosymmetric polytopes, see \cite{do05,LiNo06}. 
Decreasing $p$ while keeping $N$ fixed shows that
the  bound \eqref{Improve3s} becomes in fact irrelevant for small $p$,
since the bound $m \ge 2s$ holds as soon as there exists a pair $(A,\Delta)$ allowing exact  recovery of all $s$-sparse vectors, see \cite[Lem. 3.1]{codade09}.
Combining the two bounds, we see that $s$-sparse recovery by $\ell_p$-minimization forces
$$
m\ge C_1 s \, \big( 1+ p\ln(N/(C_2s)) \big),
$$
for some constants $C_1,C_2 > 0$.
Interestingly, if 
such an inequality is fulfilled (with possibly different constants $C_1,C_2$)  
and if $A$ is a Gaussian random matrix, then the pair $(A,\Delta_p)$ allows $s$-sparse 
recovery with high probability, see \cite{chst08}.
We note, however, that exact $\ell_p$-minimization with $p<1$, as a non-convex optimization program, encounters significant difficulties of implementation.
For more information on compressive sensing, we refer to \cite{ca06-1,carota06,cata06,codade09,do06-2,ra09-1}.

\subsection{Acknowledgements}

The first author is supported by the French National Research Agency (ANR) through the project ECHANGE (ANR-08-EMER-006).
The third and fourth author acknowledge support by the Hausdorff Center for Mathematics, University of Bonn.
The third author acknowledges funding through the WWTF project SPORTS (MA07-004). 

\section{Lower Bounds}

In this section, we use compressive sensing methods to establish the lower bound in \eqref{upper:lower:Gelfand}, hence the lower bound in \eqref{upper:lower:Gelfand2} as a by-product. 
Precisely, we show the following result,
in which the restriction $q \le 2$ is not imposed.

\begin{proposition}
\label{PropLB}
For $0<p\leq1$ and $p < q\leq \infty$,
there exists a constant $c_{p,q}>0$  such that 
\begin{equation}\label{main}
  d^m(B_p^N,\ell_q^N) \geq c_{p,q}\min\Big\{1,\frac{\ln(eN/m)}{m}\Big\}^{1/p-1/q}\quad,\quad m<N\,.
\end{equation}
\end{proposition} 
\noindent

The proof of Proposition \ref{PropLB} involves several auxiliary steps.
We start with a result \cite{grni03,grni07} on the 
unique recovery of sparse vectors 
via $\ell_p$-minimization for $0<p\leq 1$.
A proof  is included for the reader's convenience. 
We point out that,
given a subset $S$ of $[N] := \{1,...,N\}$ and a vector $v \in \R^N$,
we denote by $v_{S}$ the
vector that coincides with $v$ on $S$ and that
vanishes on the complementary set $S^c := [N]\setminus S$.
\begin{lemma}\label{NSP} Let $0<p\leq 1$ and $N,m,s \in \N$ with $m,s< N$. For a matrix
$A \in \R^{m\times N}$, the following statements are equivalent.
\begin{itemize}
\item[(a)] Every $s$-sparse vector $x$ is the unique minimizer of $\|z\|_p$ subject to
$Az = Ax$, 

\item[(b)] A satisfies the $p$-null space property of order $s$, i.e., 
for every $v\in \ker A\setminus \{0\}$ and every $S \subset [N]$ with $|S|\leq s$,  
$$
  \|v_S\|_p^p < \frac{1}{2}\|v\|_p^p\,.
$$
\end{itemize}
\end{lemma}

\begin{Proof}
$(a) \Rightarrow (b)$: Let $v\in \ker A\setminus\{0\}$ and $S \subset [N]$ with $|S| \leq s$.
Since $v = v_{S}+v_{S^{c}}$ satisfies $A v =0$,
we have $Av_{S} = A(-v_{S^c})$.
Then, since $v_{S}$ is $s$-sparse, (a) implies
$$
  \|v_S\|_p^p < \|-v_{S^{c}}\|_p^p = \|v_{S^{c}}\|_p^p.
$$
Adding $\|v_{S}\|_p^p$ on both sides and using $\|v_{S^c}\|_p^p + \|v_S\|_p^p = \|v\|_p^p$ gives (b).

$(b) \Rightarrow (a)$: Let $x$ be an $s$-sparse vector and let $S:=\supp x$. Let further $z\neq x$ be such that $Az = Ax$. 
Then $v := x-z \in \ker A\setminus \{0\}$ and 
\begin{equation}\label{eq11}
\|x\|_p^p \leq \|x_S-z_S\|_p^p + \|z_S\|_p^p = \|v_S\|_p^p + \|z_S\|_p^p\,,
\end{equation}
where the first estimate is a consequence of the $p$-triangle inequality \eqref{ptriangle}.
Clearly, (b) implies $\|v_S\|_p^p < \|v_{S^c}\|_p^p$. Plugging this into (\ref{eq11}) and using that 
$v_{S^c} = -z_{S^c}$ gives 
$$
  \|x\|_p^p < \|v_{S^c}\|_p^p + \|z_S\|_p^p = \|z_{S^c}\|_p^p + \|z_S\|_p^p = \|z\|_p^p.
$$
This proves that $x$ is the unique minimizer of $\|z\|_p$ subject to
$Az = Ax$.
\end{Proof}

\noindent
The next auxiliary step is a well-known combinatorial lemma, see for instance 
\cite{avno94,bumirave00,grsl80}, \cite[Lem.\ 3.6]{MePaRu05}.
A proof that provides explicit constants is again included for the reader's convenience.
\begin{lemma} \label{Pajor}
Let $N,s\in \N$ with $s < N$\,. There exists a family $\mathcal{U}$ of subsets  of $[N]$
such that 
\begin{itemize} 
\item[(i)] Every set in $\mathcal{U}$ consists of exactly $s$ elements.
\item[(ii)]  For all $I, J\in \mathcal{U}$ with $I\neq J$, it holds $|I\cap J| < s/2$.
\item[(iii)] The family $\mathcal{U}$ is ``large'' in the sense that 
$$
  |\mathcal{U}| \geq \Big(\frac{N}{4s}\Big)^{ s/2 }\,.
$$
\end{itemize}
\end{lemma}

\begin{Proof}
We may assume that $s \le N/4$, for otherwise we can take a family $\mathcal{U}$ consisting of just one element.
Let us denote by $\mathcal{B}(N,s)$ the family 
of subsets of $[N]$ having exactly $s$ elements. 
This family has size $|\mathcal{B}(N,s)| = \binom{N}{s}$.
We draw an arbitrary element $I_1 \in 
\mathcal{B}(N,s)$ and collect in a family $\mathcal{A}_1$ all the sets 
$J \in \mathcal{B}(N,s)$ such that $|I\cap J| \ge s/2$\,.
Then $\mathcal{A}_1$ has size at most 
\begin{equation}\label{eq10}
  \sum\limits_{k=\lceil s/2 \rceil}^{s} \binom{s}{k}\binom{N-s}{s-k} 
\le 2^s\max\limits_{\lceil s/2 \rceil\leq k\leq s \; \; }\binom{N-s}{s-k}
  = 2^s\binom{N-s}{\lfloor s/2 \rfloor},
\end{equation}
the latter inequality holding because $\lfloor s/2 \rfloor \le (N-s)/2$ when $s \le N/2$.
We throw away $\mathcal{A}_1$ and observe that  every element in $J \in \mathcal{B}(N,s)\setminus \mathcal{A}_1$ satisfies $|I_1 \cap J| < s/2$\,.
Next we draw an arbitrary element $I_2 \in \mathcal{B}(N,s)\setminus \mathcal{A}_1$, provided that the latter is not empty. 
We repeat the procedure, i.e., we define a family $\mathcal{A}_2$ relative to $I_2$ and draw an arbitrary element  $I_3 \in \mathcal{B}(N,s)\setminus(\mathcal{A}_1 \cup \mathcal{A}_2)$, and so forth until no more elements are left.
The size of each set  $\mathcal{A}_i$ can always be estimated from above by (\ref{eq10}).
This results in a collection  $\mathcal{U} = \{I_1,\hdots,I_L\}$ of subsets of $[N]$ satisfying (i) and (ii).
We finally observe that
\begin{eqnarray*}
  L & \geq & \frac{\binom{N}{s}}{2^s\binom{N-s}{\lfloor s/2 \rfloor}}
  = \frac{1}{2^s}\, \frac{N(N-1)\cdots(N-s+1)}{(N-s)(N-s-1)\cdots(N-s-\lfloor s/2 \rfloor +1)} \, \frac{1}{s(s-1)\cdots(\lfloor s/2 \rfloor + 1)}\\
  & \ge & 
  \frac{1}{2^s} \, \frac{N(N-1)\cdots (N-\lceil s/2 \rceil +1)}{s(s-1) \cdots (s-\lceil s/2 \rceil +1)}
  \ge \frac{1}{2^s} \, \Big(\frac{N}{s}\Big)^{\lceil s/2 \rceil} \ge \Big(\frac{N}{4s}\Big)^{s/2 }.
\end{eqnarray*}
This concludes the proof by establishing (iii).
\end{Proof}

\noindent
We now use Lemma \ref{Pajor} for the final auxiliary result,
which is quite interesting on its own.
It gives an estimate of the minimal number of measurements for exact recovery of sparse vectors via $\ell_p$-minimization,
where $0 < p \le 1$.

\begin{lemma}\label{pajor2} 
Let $0<p\leq 1$ and  $N,m,s \in \N$ with $m < N$ and $s < N/2$. 
If $A \in \R^{m\times N}$ is a matrix such that every $2s$-sparse vector $x$ is a minimizer of $\|z\|_p$ subject to $Az = Ax$,
then
$$
  m\geq c_1ps\ln\Big(\frac{N}{c_2s}\Big)\,,
$$
where $c_1 := 1/\ln 9 \approx 0.455$ and $c_2 := 4$\,. 
\end{lemma}

\begin{remark} Lemma 2.6 could be rephrased  (with modified constants)
by replacing $2s$-sparse vectors, $s \ge 1$, by $s$-sparse vectors, $s \ge 2$.
In the case $s=1$, it is possible for every $1$-sparse vector $x$ to be a (nonunique) minimizer of $\|z\|_1$ subject to $A z = A x$,
yet  $m \ge c_1 p \ln(N/c_2)$ fails for all constants $c_1,c_2>0$.
This can be verified by taking $m=1$ and  $A = \begin{bmatrix} 1 & 1 & \cdots & 1 & 1 \end{bmatrix}$.\end{remark}

\begin{Proof} 
We consider the quotient space 
$$
 X:=\R^N/\ker A = \{[x]:= x+\ker A\,,\,x\in \R^N\},
$$
which has algebraic dimension $r := \mbox{rank}\,A \leq m$. 
It is a quasi-Banach space equipped with 
$$
  \|[x]\|_{A,p} := \inf\limits_{v\in \ker A} \|x+v\|_p.
$$
Indeed,
a simple computation reveals that $\|\cdot\|_{A,p}$ satisfies the $p$-triangle inequality, i.e., 
$$
   \|[x]+[y]\|_{A,p}^{p} \leq \|[x]\|_{A,p}^p + \|[y]\|_{A,p}^p\,.
$$
By assumption, 
the quotient map $[\cdot]$ preserves the norm of every $2s$-sparse vector.
We now choose a family $\mathcal{U}$ of subsets of $[N]$  
satisfying (i), (ii), (iii) of Lemma \ref{Pajor}.
For a set $I \in \mathcal{U}$, we define an element $x_I \in \ell_p^N$  with $\|x_I\|_p=1$ by
\begin{equation}
\label{DefxI} 
 x_I :=  \frac{1}{s^{1/p}}\sum\limits_{i\in I}e_i\,,    
\end{equation}
where $(e_1,\ldots,e_N) $ denotes the canonical basis of $\R^N$. 
For $I,J \in \mathcal{U}, \; I \not= J$,  (ii) yields
$$
  \|x_I-x_J\|_p^p > \frac{2s-2s/2}{s} = 1\,.
$$
Since  the vector $x_I-x_J$ is a $2s$-sparse vector,
we obtain
$$
   \|[x_I]-[x_J]\|_{A,p} = \|[x_I-x_J]\|_{A,p} 
   = \|x_I-x_J\|_p > 1\,.
$$
The $p$-triangle inequality implies that
$\{[x_I] + (1/2)^{1/p}B_X,~ I\in \mathcal{U}\}$ is a disjoint collection of balls included in the ball $(3/2)^{1/p} B_X$,
where $B_X$ denotes the unit ball of $(X, \|\cdot\|_{A,p})$.
Let $\mbox{vol}(\cdot)$ denote a volume 
form on $X$, that is a translation invariant measure satisfying 
$\mbox{vol}(B_X)>0$ and $\mbox{vol}(\lambda B_X) = \lambda^r\mbox{vol}(\lambda B_X)$ for all 
$\lambda>0$ (such a measure exists since $X$ is isomorphic to $\R^r$).
The volumes satisfy the relation 
$$
\sum_{I \in \mathcal{U}} \mbox{vol} \big([x_I] + (1/2)^{1/p}B_X \big) 
\le 
\mbox{vol} \big((3/2)^{1/p} B_X \big).
$$
By translation invariance and homogeneity, we then derive 
$$
|\mathcal{U}| \, (1/2)^{r/p} \, \mbox{vol} \big(B_X \big) \leq (3/2)^{r/p} \, \mbox{vol} \big(B_X \big)\,.
$$
As a result of (iii), we finally obtain
$$
\Big( \frac{N}{4s} \Big)^{ s/2 }   \le 3^{r/p} \le 3^{m/p}.
$$
Taking the logarithm on both sides gives the desired result.
\end{Proof}

Now we are ready to prove Proposition \ref{PropLB}\,.
The underlying idea is that a small Gelfand width would imply $2s$-sparse recovery for $s$ large enough to violate the conclusion of Lemma \ref{pajor2}. \newline

\begin{Proof} 
With $c:=(1/2)^{2/p-1/q}$ and $d:= 2c_1 p / (4+c_1) \approx 0.204 \, p$, we are going to prove that
\begin{equation}
\label{LBwithCsts}
d^m(B_p^N,\ell_q^N) \ge  c \, \mu^{1/p-1/q},
\quad \mbox{ where }
\mu:= \min\Big\{1,\frac{d\ln(eN/m)}{m}\Big\}.
\end{equation}
The desired result will follow with $c_{p,q} := c \, d^{1/p-1/q}$.
By way of contradiction, we assume that
$d^m(B_p^N,\ell_q^N) < c \, \mu^{1/p-1/q}$. 
This implies the existence of a matrix $A \in \R^{m \times N}$ such that, for all $v \in \ker A \setminus \{0\}$,
$$
\|v\|_q < c \, \mu^{1/p-1/q} \|v\|_p.
$$
For a fixed $v \in \ker A \setminus \{0\}$, 
in view of the inequalities $\|v\|_p \le N^{1/p-1/q} \|v\|_q$ 
and $ c \le (1/2)^{1/p-1/q}$,
we derive $1 < (\mu N/2)^{1/p-1/q}$, so that $1 \le 1/\mu < N/2$.
We then define $s:=\lfloor 1/\mu \rfloor \ge 1$, so that
$$
\frac{1}{2 \mu} < s \le \frac{1}{\mu}.
$$
Now, for $v \in \ker A \setminus \{0\}$ and $S \subset [N]$ with $|S| \le 2s$, we have
$$
\|v_S\|_p \le (2s)^{1/p-1/q} \|v_S\|_q
\le (2s)^{1/p-1/q} \|v\|_q < c \, (2s\mu)^{1/p-1/q} \|v\|_p \le 
\frac{1}{2^{1/p}} \|v\|_p.
$$
This shows that the $p$-null space property of order $2s$ is satisfied.
Hence, Lemmas \ref{NSP} and \ref{pajor2} imply 
\begin{equation}
\label{LBform}
m\geq c_1ps\ln\Big(\frac{N}{c_2s}\Big)\,.
\end{equation}
Besides, since the pair $(A,\Delta_p)$ allows exact recovery of all $2s$-sparse vectors, we have 
\begin{equation}
\label{LBform2}
m \ge 2 \, (2s) = c_2 s.
\end{equation}
Using \eqref{LBform2} in \eqref{LBform}, it follows that
$$
m \ge c_1ps\ln\Big(\frac{N}{m}\Big)
= c_1ps\ln\Big(\frac{eN}{m}\Big) - c_1ps
>\frac{c_1 p }{2 \mu} \ln\Big(\frac{eN}{m}\Big) - \frac{c_1}{4} m.
$$
After rearrangement, we deduce 
$$
m > \frac{2c_1p}{4+c_1} \, \frac{\ln (eN/m)}{\min \big\{ 1, d \ln(eN/m) / m \big\}}
\ge \frac{2c_1p}{4+c_1} \, \frac{\ln (eN/m)}{ d \ln(eN/m) } \,m = m.
$$
This is the desired contradiction.
\end{Proof}

\begin{remark} 
When $m$ is close to $N$,
the lower estimate \eqref{LBwithCsts} is rather poor.
In this case,
a nice and simple argument proposed to us by Vyb{\'i}ral gives the improved estimate
\begin{equation}
\label{LBVyb}
d^m(B_p^N,\ell^N_q) \geq \Big(\frac{1}{m+1}\Big)^{1/p-1/q},\qquad m<N\,.
\end{equation}
Indeed, for an arbitrary matrix $A \in \R^{m \times N}$, the kernel of $A$ and the $(m+1)$-dimensional space $\{ x \in \R^N: x_i=0 \mbox{ for all } i > m+1 \}$ have a nontrivial intersection.
We then choose a vector $v \neq 0$ in this intersection,
and \eqref{LBVyb} follows from 
the inequality  $\|v\|_p \leq (m+1)^{1/p-1/q}\|v\|_q$. 
\end{remark}

We close this section with the important observation that any measurement/reconstruction scheme that provides $\ell_1$-stability requires a number of measurements scaling at least like the sparsity times a $\log$-term. 
This may be viewed as a consequence of Propositions \ref{prop:Em} and \ref{PropLB}.
Indeed, fixing $p<1$, the inequalities \eqref{StableRec} and   \eqref{compressible} 
imply
$$
d^m(B_p^N,\ell_1^N) \le E_m(B_p^N,\ell_1^N) 
\le C \sup_{x \in B_p^N} \sigma_s(x)_1 \le \frac{C }{s^{1/p-1}}.
$$
The lower bound \eqref{main} for the Gelfand width then yields, for some constant $c$,
$$
c \, \min\Big\{1,\frac{\ln(eN/m)}{m}\Big\} \le \frac{1}{s}.
$$
We derive either $s \le 1/c$ or $m \ge c s \ln(eN/m)$.
In short, if $s > 1/c$, then $\ell_1$-stability implies $m \ge c s \ln(eN/m)$
--- which can be shown to imply in turn $m \ge c' s \ln(eN/s)$.
We provide below a direct argument that removes the restriction $s > 1/c$.
It uses Lemma \ref{Pajor} and works also for $\ell_p$-stability with $p < 1$.
It borrows ideas from a paper by Do Ba et al. \cite[Thm. 3.1]{doinprwo10},
which contains the case $p=1$ in a stronger non-uniform version.

\begin{Theorem}\label{stable}
Let $N,m,s \in \N$ with $m,s<N$.
Suppose that a measurement matrix $A \in \R^{m \times N}$ and a reconstruction map $\Delta: \R^N \to \R^m$ are stable in the sense that, for all $x \in \R^N$,
$$
\|x-\Delta(Ax)\|_p^p \le C \sigma_s(x)_p^p
$$
for some constant $C > 0$ and some $0<p\leq 1$.
Then there exists a constant $C'>0$ depending only on $C$  such that
$$
m \ge C' p \, s \ln(eN/s).
$$
\end{Theorem}

\begin{Proof} 
We consider again a family $\mathcal{U}$ of subsets of $[N]$ given by Lemma \ref{Pajor}.
For each $I \in \mathcal{U}$,
we define an $s$-sparse vector $x_I$ with $\|x_I\|_p=1$ as in \eqref{DefxI}. 
With $\rho:=(2(C+1))^{-1/p}$,
we claim that 
$\{A(x_I+ \rho B_p^N), I \in \mathcal{U}\}$ is a disjoint collection of subsets of $A(\R^N)$, which has algebraic dimension $r \le m$.
Suppose indeed that there exist $I,J \in \mathcal{U}$ with $I \not= J$
and $z,z' \in \rho B_p^N$ such that $A(x_I+z) = A(x_J+z')$.
A contradiction follows from
\begin{eqnarray*}
1  <  \|x_I-x_J\|_p^p & \le & 
\| x_I + z - \Delta(A(x_I+z))\|_p^p + \|x_J + z' - \Delta(A(x_J+z')) \|_p^p + \|z\|_p^p + \|z'\|_p^p\\
& \le & C \sigma_s(x_I+z)_p^p +  C \sigma_s(x_J+z')_p^p +   \|z\|_p^p + \|z'\|_p^p\\
& \le & C \|z\|_p^p + C \|z'\|_p^p +  \|z\|_p^p + \|z'\|_p^p \le 2(C+1) \rho^p = 1.
\end{eqnarray*}
We now observe that the collection $\{A(x_I+ \rho B_p^N), I \in \mathcal{U}\}$ is contained in $(1+\rho^p)^{1/p} A(B_p^N)$.
As in the proof of Lemma \ref{pajor2},
we use a standard volumetric argument to derive
$$
|\mathcal{U}| \, \rho^r \,  \mbox{vol}\big(A(B_p^N)\big) = 
\sum_{I \in \mathcal{U}} \mbox{vol}\big(A(x_I+ \rho B_p^N)\big)
\le \mbox{vol}\big((1+\rho^p)^{1/p} A(B_p^N)\big) 
= (1+\rho^p)^{r/p} \mbox{vol}\big(A(B_p^N)\big).
$$
We deduce that
$$
\Big( \frac{N}{4s} \Big)^{ s/2 }
\le (\rho^{-p}+1)^{r/p} \le  (\rho^{-p}+1)^{m/p} = (2C+3)^{m/p}.
$$
Taking the logarithm on both sides yields
$$
m \ge c p s \ln(N/(4s)),
\qquad \mbox{ with } c:=1/(2 \ln(2C+3)). 
$$
Finally, noticing that $m \ge 2s$ because the pair $(A,\Delta)$ allows exact $s$-sparse recovery, we obtain
$$
m \ge c p s \ln(eN/s) - c p s \ln(4e) \ge c p s \ln(eN/s)  - \frac{c \ln(4e)}{2} m.
$$
The desired result follows with $C':=(2c)/(2+c \ln(4e))$. 
\end{Proof}

\section{Upper Bounds}

In this section, we establish the upper bound in \eqref{upper:lower:Gelfand2},
hence the upper bound in \eqref{upper:lower:Gelfand} as a by-product. 
As already mentioned in the introduction,
the bound for the Gelfand width of $\ell_p$-balls was already provided by Vyb{\'i}ral in \cite{vy08}, but the bound for the Gelfand width of weak-$\ell_p$-balls  is indeed new.

\begin{proposition}
For $0<p < 1$ and $p < q\leq 2$, 
there exists a constant $C_{p,q}>0$ such that 
\begin{equation}\label{f10}
  d^m(B_{p,\infty}^N,\ell_q^N) \leq C_{p,q}\min\Big\{1,\frac{\ln(eN/m)}{m}\Big\}^{1/p-1/q}\quad,\quad m<N\,.
\end{equation}
\end{proposition}

The argument relies again on compressive sensing methods.
According to Proposition \ref{prop:Em}, 
it is enough to establish the upper bound for the quantity $E_m(B_{p,\infty}^N,\ell_q^N)$.
This is done in the following theorem, which we find rather illustrative because it shows that, even when $p < 1$, an optimal reconstruction map $\Delta$ 
for the realization of the number $E_m(B_{p,\infty}^N,\ell_q^N)$ can be chosen to be the $\ell_1$-minimization mapping, at least when $q \geq  1$. 
The argument is originally due to Donoho for $q=2$ \cite[Proof of Theorem 9]{do06-2} and can be extended to all $2\geq q >p$.

\begin{Theorem}\label{thm:optdec} For $0<p < 1$ and $p<q\leq 2$, there exists a matrix $A \in \R^{m \times N}$ such that, with $r = \min\{1,q\}$,
\[
\sup_{x \in B_{p,\infty}^N} \|x - \Delta_r(Ax)\|_q \leq C_{p,q} \min \Big\{1,\frac{\ln(N/m)+1}{m} \Big\}^{1/p-1/q},
\]
where $C_{p,q}>0$ is a constant that depends only on $p$ and $q$.
\end{Theorem}

\begin{Proof} 
Let $C_1$ be the constant in \eqref{mRIP} relative to the RIP 
associated with $\delta=1/3$, say. We choose a constant $D>0$ large enough to have
$$
D/2>e,
\qquad
\frac{D/2}{1+\ln(D/2)} > C_1.
$$
We are going to prove that, for any $x \in B_{p,\infty}^N$,
\begin{equation}
\label{LastPf}
\|x - \Delta_r(Ax)\|_q \le C'_{p,q} \min \Big\{1,\frac{D \ln(eN/m)}{m} \Big\}^{1/p-1/q}
\end{equation}
for some constant $C'_{p,q} > 0$. This will imply the desired result with $C_{p,q}:=C'_{p,q} D^{1/p-1/q}$.\\
{\em Case 1:} $\displaystyle{m > D \ln (eN/m)}$.\\
We define $s \ge 1$ as the largest integer smaller than $ m/ (D \ln (eN/m))$, so that
\begin{equation}
\label{DefS}
\frac{m}{2 D \ln(e N/ m)} \le s < \frac{m}{D \ln(e N/ m)}.
\end{equation}
Putting $t=2s$ and noticing that $t/m < 2/D < 1/e$ and that $u \mapsto u \ln(u)$ is decreasing on $ [0,1/e]$, we obtain
$$
m > \frac{D}{2} t \ln(eN/m) = \frac{D}{2} t \ln(eN/t) + \frac{D}{2} \, m\, (t/m) \ln(t/m)
>  \frac{D}{2} t  \ln(eN/t) - m \, \ln(D/2),
$$
so that
$$
m > \frac{D/2}{1+\ln(D/2)} t  \ln(eN/t)  > C_1 t  \ln(eN/t). 
$$
It is then possible to find a matrix $A \in \R^{m \times N}$ with $\delta_{t}(A) \le \delta$.
In particular, we have $\delta_{s}(A) \le \delta$.
Now, given $v := x - \Delta_r(Ax) \in \ker A$, we decompose $[N]$ as the disjoint union of sets $S_1,S_2,S_3,\ldots$ of size $s$  (except maybe the last one) in such a way that
$|v_i| \ge |v_j|$ for all $i \in S_{k-1}$, $j \in S_k$, and $k \ge 2$.
This easily implies $\big( \|v_{S_k}\|_2^2 / s \big)^{1/2} \le \big( \|v_{S_{k-1}}\|_r^r/s \big)^{1/r}$, i.e.,
\begin{equation}
\label{CompSkSk-1}
\|v_{S_k}\|_2 \le \frac{1}{s^{1/r-1/2}} \|v_{S_{k-1}}\|_r,
\qquad k \ge 2.
\end{equation}
Using the $r$-triangle inequality, we have
$$
\|v\|_q^r = \Big\| \sum_{k \ge 1} v_{S_k} \Big\|_q^r
\le \sum_{k \ge 1} \| v_{S_k}\|_q^r 
\le \sum_{k \ge 1} \big(s^{1/q-1/2} \| v_{S_k}\|_2 \big)^{r}
\le \sum_{k \ge 1} \Big(\frac{s^{1/q-1/2}}{\sqrt{1-\delta}} \| A v_{S_k}\|_2\Big)^r. 
$$
The fact that $v \in  \ker A$ implies $Av_{S_1} =  -\sum_{k \ge 2} A v_{S_k} $.
It follows that
\begin{eqnarray*}
\|v\|_q^r & \le &
\Big(\frac{s^{1/q-1/2}}{\sqrt{1-\delta}} \Big)^r \Big( \sum_{k \ge 2}  \| A v_{S_k}\|_2 \Big)^r
+ \Big(\frac{s^{1/q-1/2}}{\sqrt{1-\delta}} \Big)^r \sum_{k \ge 2}  \| A v_{S_k}\|_2^r\\
& \le & 2 \Big(\frac{s^{1/q-1/2}}{\sqrt{1-\delta}} \Big)^r \sum_{k \ge 2}  \| A v_{S_k}\|_2^r
\le 2 \Big(\sqrt{\frac{1+\delta}{1-\delta}} s^{1/q-1/2}\Big)^r  \sum_{k \ge 2} \| v_{S_k}\|_2^r.
\end{eqnarray*}
We then derive, using the inequality \eqref{CompSkSk-1}, 
$$
\|v\|_q^r \le 
2 \Big(\sqrt{\frac{1+\delta}{1-\delta}} \frac{1}{s^{1/r-1/q}}\Big)^{r} \sum_{k \ge 1} \| v_{S_k}\|_r^r .
$$
In view of the choice $\delta=1/3$ and of \eqref{DefS}, we deduce
\begin{equation}
\label{Eq1LastPf}
\|x-\Delta_r(Ax)\|_q \le 2^{1/r} \sqrt{2} \Big( \frac{2D \ln(eN/m)}{m} \Big)^{1/r-1/q} \|x-\Delta_r(Ax)\|_r\,.
\end{equation}
Moreover, in view of $\delta_{2s} \leq 1/3$ and of Theorem \ref{thm:l1}, there exists a constant $C>0$ such that 
\begin{equation}
\label{Eq2LastPf}
\|x - \Delta_r(Ax)\|_r \le C^{1/r} \sigma_s(x)_r.
\end{equation}
Finally, using \eqref{compressible2} and \eqref{DefS}, we have 
\begin{equation}
\label{Eq3LastPf}
\sigma_s(x)_r \le  \, \frac{D_{p,r}}{s^{1/p-1/r}} 
\le \, D_{p,r}\Big( \frac{2D \ln(eN/m)}{m} \Big)^{1/p-1/r}\,.
\end{equation}
Putting \eqref{Eq1LastPf}, \eqref{Eq2LastPf}, and \eqref{Eq3LastPf} together,
we obtain, for any $x \in B_{p,\infty}^N$, 
$$
\|x - \Delta_r(Ax)\|_q \le C''_{p,q}\Big( \frac{D \ln(eN/m)}{m} \Big)^{1/p-1/q}
=C''_{p,q}\min \Big\{1,\frac{D \ln(eN/m)}{m} \Big\}^{1/p-1/q},
$$
where $C''_{p,q}:=C^{1/r}D_{p,r}2^{1/r+1/2+1/p-1/q}$.\\
{\em Case 2:} $\displaystyle{m \le D \ln (eN/m )}$.\\
We simply choose the matrix $A \in \R^{m \times N}$ as $A=0$.
Then, for any $x \in B_{p,\infty}^N$, we have
$$
\|x - \Delta_r(Ax)\|_q = \|x\|_q \le C'''_{p,q} \|x\|_{p,\infty} \le C'''_{p,q},
$$
for some constant $C'''_{p,q}>0$.
This yields
$$
\|x - \Delta_r(Ax)\|_q \le  C'''_{p,q} \min \Big\{1,\frac{D \ln(eN/m)}{m} \Big\}^{1/p-1/q}.
$$
Both cases show that \eqref{LastPf} is valid with $C'_{p,q}:=\max\{C''_{p,q},C'''_{p,q}\}$.
This completes the proof.
\end{Proof}

\begin{remark}
The case $p=1$, for which $r=1$, is not covered by our arguments.
Since $\sup_{x\in B^N_{1,\infty}}\sigma_s(x)_1 \asymp \log(N/s)$ the quantity $\sigma_s(x)_1$ cannot be bounded by a constant times $\|x\|_{1,\infty}$ in order to obtain \eqref{Eq3LastPf}. Instead, the additional log-factor $\log(N/m)$ appears on the right-hand side and therefore in the upper estimate of \eqref{upper:lower:Gelfand2} in the case $p=1$. The correct behavior of the Gelfand widths of weak-$\ell_1$-balls does not seem to be known. 
Nonetheless, the inequality $\sigma_s(x)_1 \le \|x\|_1$ is always true. This yields the well-known upper estimate for the Gelfand widths of $\ell_1$-balls and hence completes the proof of Theorem \ref{thm:main}.
\end{remark}

\end{document}